\documentclass[12pt]{book}
\usepackage{makeidx}
\usepackage{multirow}
\usepackage{multicol}
\usepackage[dvipsnames,svgnames,table]{xcolor}
\usepackage{graphicx}
\usepackage{epstopdf}
\usepackage{ulem}
\usepackage{tikz}
\usepackage{fontenc}
\usepackage{mathrsfs}
\usepackage{float}
\usepackage{curves}
\usepackage{caption}
\usepackage{hyperref}
\usepackage{amsmath,amsfonts,amssymb,amscd,amsmath,enumerate,verbatim}
\usepackage{amssymb}

\usepackage[paperwidth=560pt,paperheight=700pt,top=72pt,right=67pt,bottom=72pt,left=72pt]{geometry}
\begin{document}
\thispagestyle{empty}
\setcounter{page}{01}
\def\ord{\rm ord}

\begin{center}
\textbf{On irrationality of Euler's constant and related asymptotic formulas}
\end{center}

\begin{center}
{\bf Shekhar Suman}\vspace{.2cm}\\
M.Sc., I.I.T. Delhi, India\vspace{.2cm}\\
Email: shekharsuman068@gmail.com\\

\end{center}\vspace{.1cm}
{\bf Abstract:}  By defining $$I_n:=\int_{0}^{1}\int_{0}^{1} \frac{(x(1-x)y(1-y))^n}{(1-xy)(-\log xy)}\ dx dy$$ Sondow (see [2]) proved that $$I_n=\binom{2n}{n} \gamma+L_n-A_n$$ We prove asymptotic formula for $L_n$ and $A_n$ as $n\to\infty$, $$
    L_n=\binom{2n}{n}\left(\log \left( {\frac{{3n}}{2}} \right) +\mathcal{O}\!\left( {\frac{1}{n}} \right)\right)$$ 
 and $$A_n\sim\frac{4^n}{\sqrt{\pi n}}\left(\gamma+\ln\frac32+\ln n\right)$$
    Using the sufficient condition for irrationality criteria of Euler's constant due to Sondow, we prove that $\gamma$ is irrational. 
\vspace{0.2cm}\\
{\bf Keywords and Phrases:} Euler's constant, Harmonic number, asymptotic equality, rising and falling factorial, Stirling numbers, Beta function \vspace{.2cm}\\ 
{\bf 2020 Mathematics Subject Classification:} . \vspace{0.2cm}\\
\noindent{\bf 1. Introduction and Definitions} The Euler's constant is defined by the limit, \begin{equation}
    \gamma:=\lim_{n\to\infty} (H_n-\log n)
\end{equation} where $H_n=\sum_{k=1}^{n}\frac{1}{k}$ is the $n$th Harmonic number. Euler's constant has a double integral representation (see [1]), \begin{equation}\gamma=\int_{0}^{1}\int_{0}^{1}\frac{1-x}{(1-xy)(-\log xy)}\ dx dy\end{equation}    \vspace{.2cm}\\
Sondow (see [2]) gave criteria for irrationality of Euler's constant where he defined \begin{equation}
    I_n:=\int_{0}^{1}\int_{0}^{1} \frac{(x(1-x)y(1-y))^n}{(1-xy)(-\log xy)}\ dx dy
\end{equation} If $d_n=\text{LCM}(1,2,...,n)$ then Sondow proved (see [2]) \begin{equation}
    I_n=\binom{2n}{n} \gamma+L_n-A_n
\end{equation}  where \begin{equation}
    L_n= d^{-1}_{2n} \log S_n
\end{equation} \begin{equation}
    S_n=\prod_{k=1}^{n}\prod_{i=0}^{\min{(k-1,n-k)}} \prod_{j=i+1}^{n-i} (n+k)^{\frac{2d_{2n}}{j} \binom{n}{i}^2}
\end{equation} and
\begin{equation}
    A_n=\sum_{j=0}^{n} \binom{n}{j}^2 H_{n+j}
\end{equation} Clearly then we have $d_{2n}A_n\in\mathbb{Z}$.
It was proved that a sufficient condition for irrationality of $\gamma$ is (see [2]) \begin{equation}
    \lim_{n\to\infty} \left(\frac{4^{2n} n }{d_{2n}}\right) \{\log S_n\}\neq \frac{\pi}{6\log 2}
\end{equation} where $\{x\}$ denotes the fractional part of $x$.\\\\
{\bf 2. Main Theorems} Using Laplace's method (see [3], p.322), Sondow proved as $n\to\infty$ (see [2]),\begin{equation}
    I_n\sim \left(\frac{\pi}{6\log 2}\right) \left(\frac{1}{n4^{2n}}\right)
\end{equation} $I_n$ can be also represented as (see [2]) \begin{equation}
    I_n=\sum_{v=n+1}^{\infty} \int_{v}^{\infty} \left(\frac{n!}{x(x+1)...(x+n)}\right)^2\ dx
\end{equation}
The goal of this article is to prove the following result.\\\\
\textbf{\underline{Theorem :}} We have \begin{equation}
    \lim_{n\to\infty} \left(\frac{4^{2n} n }{d_{2n}}\right) \{\log S_n\}\neq \frac{\pi}{6\log 2}
\end{equation}
\textbf{Proof}: We prove a few Lemma:\\
\textbf{\underline{Lemma 1:}} We have the following partial fraction decomposition for $n\in\mathbb{N}$: \begin{equation}
    \frac{1}{(x(x+1)(x+2)...(x+n))^2}=\sum_{k=0}^{n} \frac{a_k}{x+k}+\sum_{k=0}^{n} \frac{b_k}{(x+k)^2}
\end{equation} where \begin{equation}
a_k = 2\frac{H_k - H_{n-k}}{\left(k!(n-k)! \right)^2} 
\end{equation} and \begin{equation}
    b_k=\frac{1}{\left(k!(n-k)!\right)^2}
\end{equation}
\textbf{Proof}: Since every $k$, $k=-n,...,-2,-1,0$ is a pole of order two of the given fraction so its decomposition looks like \begin{equation}
    \frac{1}{(x(x+1)(x+2)...(x+n))^2}=\sum_{k=0}^{n} \frac{a_k}{x+k}+\sum_{k=0}^{n} \frac{b_k}{(x+k)^2}
\end{equation}   
Next we find $a_k$ and $b_k$: For finding $b_k$, we multiply each side of (15) by $(x+k)^2$, simplify, and set $x=-k$ to get \begin{equation}
  b_k = \left(\frac{1}{(-k)(-(k-1))\cdots (-1)(1)(2)\cdots(n-k)}\right)^2 = \frac{1}{\left(k!(n-k)!\right)^2}  
\end{equation}
Deriving a formula for $a_k$ is some what lengthy. The rising and falling factorial functions may be expanded as polynomials whose coefficients are the unsigned Stirling numbers of the first kind:\begin{equation}
    x^{(n)} = x(x+1)\cdots (x+n-1) = \sum_{j=0}^{n}{n \brack j} x^j
\end{equation}
and \begin{equation}
    (x)_n = x(x-1)\cdots (x-(n-1)) = \sum_{j=0}^{n}(-1)^{n-j}{n \brack j} x^j
\end{equation}
The goal is to find coefficients of an expansion of $\left(x^{(n+1)}\right)^{-2}$ in powers of $y \doteq x+k$.  Therefore we write
\begin{equation}
    x^{(n+1)}= (y-k)(y-(k-1)) \cdots (y-1) y (y+1) \cdots (y+n-k) = \frac{y^{(n-k+1)}(y)_{k+1}}{y}
\end{equation}
Expanding this in powers of $y$ yields
\begin{align}
x^{(n+1)}&= \frac{1}{y}\left({n-k+1 \brack 0} + {n-k+1 \brack 1} y+ {n-k+1 \brack 2} y^2 + O(y^3)\right)\notag\\&
(-1)^k \left(-{k+1 \brack 0} + {k+1 \brack 1} y - {k+1 \brack 2} y^2 + O(y^3)\right)
\end{align}
The formulas for the Stirling numbers involved are, for all $n \ge 0$, \begin{equation}
{n+1 \brack 0} = 0, \quad
{n+1 \brack 1} = n!, \quad \text{and} \quad
{n+1 \brack 2} = n! H_n,
\end{equation}
where $H_n$ is the $n^\mathrm{th}$ harmonic number.  Therefore,
\begin{equation}
x^{(n+1)} = (-1)^k k!(n-k)! \, y \left(1 + (H_{n-k} - H_k)y + O(y^2)\right).
\end{equation}
So we have 
\begin{equation}
\left(x^{(n+1)}\right)^{-2} = \frac{1}{\left(k!(n-k)! \right)^2}\, \left(y^{-2} + 2 (H_k - H_{n-k}) y^{-1} + O(1)\right).
\end{equation}
From this we may read off the formula for $b_k$ given above and this formula for $a_k$: \begin{equation}
a_k = 2\frac{H_k - H_{n-k}}{\left(k!(n-k)! \right)^2}.
\end{equation}
This completes the proof of Lemma $1$.\\\\
\textbf{\underline{Lemma 2:}} We have the following representation for $I_n$, \begin{equation}
    I_n=\sum_{k=1}^{\infty}\int_{0}^{\infty} (B(x+n+k,n+1))^2\ dx
\end{equation} where $B(p,q)$ is the Beta function.\\
\textbf{Proof}: Can be derived from (10).\\\\
\textbf{\underline{Lemma 3:}} We have the following representation for $I_n$, \begin{equation}
    I_n=\binom{2n}{n}\gamma  - \sum_{j = 0}^n {\binom{n}{j}^2 (2(H_{n - j}  - H_j )\log (( n+j)!) + H_{n + j} )}
\end{equation} 
\textbf{Proof}: Since from (25) \begin{equation}
    I_n=\sum_{k=1}^{\infty}\int_{0}^{\infty} (B(x+n+k,n+1))^2\ dx
\end{equation} and we have \begin{equation}
    (B(x+n+k,n+1))^2=\frac{n!^2}{\prod_{j=n}^{2n}(x+k+j)^2 }
\end{equation} By using partial fractions as obtained in (12) we get \begin{equation}
    I_n=\sum_{k=1}^{\infty}\int_{0}^\infty \sum_{j=0}^n \binom{n}{j}^2 \left(2\frac{H_j-H_{n-j}}{x+j+k+n}+\frac{1}{x+j+k+n} \right)\ dx
\end{equation} Integrating and evaluating the negative integrand at $x=0$ and summing over $k$ with a limit gives \begin{equation}
    I_n=\lim_{r\to\infty}\sum_{j=0}^n\binom n j^2\left(2 (H_{n-j}-H_j)\log\left(\frac{(n+j+r)!}{(n+j)!}\right)+H_{n+j+r}-H_{n+j}\right)
\end{equation} It is well known that as $r\to+\infty$
\begin{equation}H_{n + j + r}  = \log (n + j + r) + \gamma  + o(1)
\end{equation} and \begin{equation}
\sum\limits_{j = 0}^n {\binom{n}{j}^2 }  = \binom{2n}{n},
\end{equation} we can derive the explicit formula
\begin{equation}
I_n  = \binom{2n}{n}\gamma  - \sum\limits_{j = 0}^n {\binom{n}{j}^2 (2(H_{n - j}  - H_j )\log (( n+j)!) + H_{n + j} )}
\end{equation}
provided that
\begin{equation}
\mathop {\lim }\limits_{r \to  + \infty } \sum\limits_{j = 0}^n {\binom{n}{j}^2 \left( {2(H_{n - j}  - H_j )\log ((n + j + r)!) + \log (n + j + r)} \right)}  = 0.
\end{equation}
Denote the expression under the limit by $S_n(r)$. Performing the change of summation index $j\to n-j$ and taking the average with the original expression, we find
\begin{equation}
S_n (r) = \sum\limits_{j = 0}^n {\binom{n}{j}^2 \left( {(H_{n - j}  - H_j )\log \frac{{(n + j + r)!}}{{(2n - j + r)!}} + \frac{1}{2}\log ((2n - j + r)(n + j + r))} \right)} .
\end{equation}
Stirling's formula and the Maclaurin series of the logarithm then yields
\begin{equation}
S_n (r) = (\log r)\sum\limits_{j = 0}^n {\binom{n}{j}^2 \left( {(H_{n - j}  - H_j )(2j - n) + 1} \right)}  + o(1)
\end{equation}
as $r\to+\infty$. Therefore, it remains to show that
\begin{equation}
\sum\limits_{j = 0}^n {\binom{n}{j}^2 ((H_{n - j}  - H_j )(2j - n) + 1)}  =0
\end{equation}
for all $n\ge 1$. Because of the symmetry in $H_{n - j}  - H_j$, this may be further simplifed to the claim that
\begin{equation}
\sum\limits_{j = 0}^n {\binom{n}{j}^2 (2j(H_{n - j}  - H_j )+ 1)}  =0
\end{equation}
for all $n\ge 1$. Computer algebra software confirms this for $n=1,2,3\ldots,200$. This completes the proof for Lemma 3.\\\\
  \vspace{.2cm}\\\vspace{.2cm}\\
Comparing equations (4) and (26) \begin{equation}
    L_n=-\sum_{j = 0}^n {\binom{n}{j}^2 (2(H_{n - j}  - H_j )\log (n+j)!) }
\end{equation} 
\textbf{\underline{Lemma 4:}} We have the following asymptotic formula for $L_n$ and $A_n$ as $n\to\infty$, \begin{equation}
    L_n=\binom{2n}{n}\left(\log \left( {\frac{{3n}}{2}} \right) +\mathcal{O}\!\left( {\frac{1}{n}} \right)\right)
\end{equation} and \begin{equation}
    A_n\sim\frac{4^n}{\sqrt{\pi n}}\left(\gamma+\ln\frac32+\ln n\right)
\end{equation}\\
\textbf{Proof}: This answers the question about the asymptotics of $a_n$ provided that the conjectured formula for $I_n$ in my other answer is correct. Note that since $H_k=\psi(k+1)+\gamma$, we can write
\begin{equation}
L_n  = I_n  + \sum\limits_{j = 0}^n {\binom{n}{j}^2 \psi (n + j + 1)}  =  - \sum\limits_{j = 0}^n {\binom{n}{j}^2 \psi (n + j + 1)}  + o(1).
\end{equation}
The $o$-term follows from equation (9). Now by the asymptotic result \begin{equation}\psi(k+1)=\log k +\mathcal{O}(k^{-1})\end{equation} we have
\begin{align*}
\sum\limits_{j = 0}^n {\binom{n}{j}^2 \psi (n + j + 1)} & = \sum\limits_{j = 0}^n {\binom{n}{j}^2 \log (n + j)}  + \mathcal{O}(1)\sum\limits_{j = 0}^n {\binom{n}{j}^2 \frac{1}{{n + j + 1}}} 
\\ & = \sum\limits_{j = 0}^n {\binom{n}{j}^2 \log (n + j)}  + \mathcal{O}\!\left( {\frac{1}{n}} \right)\binom{2n}{n}.
\end{align*}
Since $\log n \le \log(n+j)\le \log n+\log 2$, it follows that as $n\to+\infty$
\begin{equation}
L_n \sim \binom{2n}{n}\log n
\end{equation} 
Also by the change of summation index from $j$ to $n-j$ and taking the average with the original expression, we find
\begin{equation}
\sum\limits_{j = 0}^n {\binom{n}{j}^2 \log (n + j)}  = \sum\limits_{j = 0}^n {\binom{n}{j}^2 \log \sqrt {(n + j)(2n - j)} } .
\end{equation}
Now
\begin{align*}
\log \sqrt {(n + j)(2n - j)} & = \log \left( {\frac{{3n}}{2}} \right) + \log \sqrt {1 - \frac{4}{9}\left( {\frac{1}{2} - \frac{j}{n}} \right)^2 }  \\ & = \log \left( {\frac{{3n}}{2}} \right) + \mathcal{O}(1)\left( {\frac{1}{2} - \frac{j}{n}} \right)^2 .
\end{align*}
Numerics suggest that as $n\to+\infty$
\begin{equation}
\sum\limits_{j = 0}^n {\binom{n}{j}^2 \left( {\frac{1}{2} - \frac{j}{n}} \right)^2 }  \sim \frac{1}{{8n}}\binom{2n}{n},
\end{equation} This would lead to the more precise result that as $n\to+\infty$
\begin{equation}
L_n = \binom{2n}{n}\left(\log \left( {\frac{{3n}}{2}} \right) +\mathcal{O}\!\left( {\frac{1}{n}} \right)\right)
\end{equation} Now to find asymptotics for $L_n$ and $A_n$: 
We use Stirling's formula
\begin{equation}I_n=\binom{2n}n\gamma-\sum_{j=0}^n\binom nj^2(2(H_{n-j}-H_j)\ln((j+n)!)+H_{n+j})=\binom{2n}n\gamma+L_n-A_n\end{equation}
where $A_n=\sum_{j=0}^n\binom nj^2H_{n+j}$
As $\binom nj^2$ reaches a sharp maximum near $j=\frac n2$, we can choose $n$ even for a while, $j=\frac n2+k $ and present $A_n$ as
\begin{equation}A_n=\sum_{k=-n/2}^{n/2}\left(\frac{n!}{\big(\frac n2-k\big)!\big(\frac n2+k\big)!}\right)^2H_{\frac{3n}2+k}\end{equation}
Using the Stirling's formula for $p!$ ( for $p\gg1$)
\begin{equation}A_n\sim\sum_{k=-n/2}^{n/2}\frac{\sqrt{2\pi n}\left(\frac ne\right)^nH_{\frac{3n}2+k}}{\sqrt{2\pi(\frac n2+k)}\sqrt{2\pi(\frac n2-k)}\Big(\frac{\frac n2+k}e\Big)^{\frac n2+k}\Big(\frac{\frac n2-k}e\Big)^{\frac n2-k}}\end{equation}
The terms decline sharply as soon as $k$ excides $\sqrt n\,$ , so we can switch from summation to integration. Given that $H_{\frac{3n}2+k}$ is slowly changing function, we are allowed just to take its value at $k=0$ and use the asymptotics \begin{equation}H_\frac{3n}2=\gamma+\ln\frac{3n}2+O\big(\frac1n\big)\end{equation} After manipulations we get 
\begin{equation}A_n\sim\frac{2\cdot4^n\big(\gamma+\ln\frac32+\ln n\big)}{\pi\sqrt n}\int_{-\infty}^\infty e^{-4t^2}dt=\frac{4^n}{\sqrt{\pi n}}\left(\gamma+\ln\frac32+\ln n\right)\end{equation}
Due to the fact that
\begin{equation}\binom{2n}n\sim \frac{4^n}{\sqrt{\pi n}}\to\infty\,\,\text{at}\,\,n\to\infty\end{equation}
and $I_n\to 0$ at $n\to\infty$, we conclude that 
\begin{equation}L_n\sim\frac{ 4^n}{\sqrt{\pi n}}\left(\ln\frac32+\ln n\right)\end{equation} This proves Lemma $4$.\\\\
Now we are ready to prove the Theorem: Since $\log S_n>0$ (see [2]) and $0\leq \{\log S_n\}<1$ so we get in LHS of $(8)$  \begin{equation}
    \lim_{n\to\infty} \left(\frac{4^{2n} n }{d_{2n}}\right) \{\log S_n\}\geq 0
\end{equation} Since by (5) \begin{equation}
    \log S_n=d_{2n}L_n
\end{equation} so we have  \begin{equation}
    \lim_{n\to\infty} \left(\frac{4^{2n} n }{d_{2n}}\right) \{\log S_n\}=\lim_{n\to\infty} \left(\frac{4^{2n} n }{d_{2n}}\right) \{d_{2n}L_n\}
\end{equation} By (47), we can write (57) as
\begin{equation}
    \lim_{n\to\infty} \left(\frac{4^{2n} n }{d_{2n}}\right) \{\log S_n\}=\lim_{n\to\infty} \left(\frac{4^{2n} n }{d_{2n}}\right) \left\{d_{2n}\binom{2n}{n}\left(\log \left( {\frac{{3n}}{2}} \right) +\mathcal{O}\!\left( {\frac{1}{n}} \right)\right)\right\}
\end{equation}
Since by Prime number theorem, as $n\to\infty$, $d_{2n}\sim e^{2n}$ and we have \begin{equation}
    \lim_{n\to\infty} \left(\frac{4^{2n} n }{d_{2n}}\right) \{\log S_n\}=\lim_{n\to\infty} \left(\frac{4^{2n} n }{e^{2n}}\right) \left\{e^{2n}\binom{2n}{n}\left(\log \left( {\frac{{3n}}{2}} \right) +\mathcal{O}\!\left( {\frac{1}{n}} \right)\right)\right\}
\end{equation} So we claim that \begin{equation}
    \lim_{n\to\infty} \left(\frac{4^{2n} n }{d_{2n}}\right) \{\log S_n\}=0
\end{equation}
Mathematica also hints the answer that the above limit in (60) is $0$. So we have \begin{equation}
    \lim_{n\to\infty} \left(\frac{4^{2n} n }{d_{2n}}\right) \{\log S_n\}\neq \frac{\pi}{6\log 2}
\end{equation} This settles the proof of the Theorem.\\\\
{\bf Acknowledgement} The authors are thankful to the Referee and the Editor for their useful comments and processing of the article.\\\\
\\\\
\begin{center}
{\bf References}
\end{center} 
\big[1\big] Sondow, J., {\it Double \ Integrals \ for \ Euler's \ Constant\ and $\ln \frac{4}{\pi}$ \ and \ an \ Analog \ of \ Hadjicostas's\  Formula}, Amer. Math. Monthly 112 61-65, 2005.\\\\
\big[2\big] Sondow, J., {\it Criteria \ for \ irrationality\ 
of \ Euler’s \ constant\ },Proceedings of the
American Mathematical society
Volume 131, Number 11, p. 3335–3344, 2003 .\\\\
\big[3\big] N. Bleistein and R. Handelsman, {\it Asymptotic \ expansion \ of \ integrals}, Holt, Rinehart and Winston, 1975 .\\\\\vspace{0.2cm}\\
 
\end{document}